\documentstyle[12pt,amssymb]{amsart}
\title[Unimodular invariants of tori]
{Unimodular invariants of totally real tori in ${\Bbb C}^n$}
\thanks{Partially supported by NSF grant DMS-9304580 at the Institute for Advanced Study}
\author{Xianghong  Gong}
\address{Institute for Advanced Study, School of  Mathematics, 
Princeton, NJ 08540}
\email{gong@@math.ias.edu}

\newtheorem{thm}{Theorem}[section]
\newtheorem{cor}[thm]{Corollary}
\newtheorem{prop}[thm]{Proposition}
\newtheorem{define}[thm]{Definition}
\newtheorem{lemma}[thm]{Lemma}

\newcommand{\be}[1]{\begin{equation}\label{eq:#1}}
\newcommand{\nee}{\end{equation}}
\newcommand{\ben}{\begin{equation*}}
\newcommand{\een}{\end{equation*}}
\newcommand{\nbeq}{\begin{eqnarray}}
\newcommand{\eeq}{\end{eqnarray}}
\newcommand{\bc}[1]{\begin{cor}\label{cor:#1}}
\newcommand{\ec}{\end{cor}}
\newcommand{\bt}[1]{\begin{thm}\label{thm:#1}}
\newcommand{\et}{\end{thm}}
\newcommand{\bl}[1]{\begin{lemma}\label{lemma:#1}}
\newcommand{\el}{\end{lemma}}
\newcommand{\bd}[1]{\begin{define}\label{lemma:#1}}
\newcommand{\ed}{\end{define}}
\newcommand{\bp}[1]{\begin{prop}\label{prop:#1}}
\newcommand{\ep}{\end{prop}}

\newcommand{\rl}[1]{Lemma~\ref{lemma:#1}}

\newcommand{\re}[1]{(\ref{eq:#1})}
\newcommand{\rp}[1]{Proposition~\ref{prop:#1}}
\newcommand{\rt}[1]{Theorem~\ref{thm:#1}}

\begin{document}
\begin{abstract}
We study the global invariants of real analytic manifolds in
the complex space with respect to the group of holomorphic unimodular 
transformations. We consider only totally real manifolds which admits a 
certain fibration over the circle. We find a complete set of invariants for
totally real tori in ${\Bbb C}^n$ which are close to the standard torus. 
The invariants are obtained by an analogous classification of complex-valued 
analytic $n$-forms on the standard torus. We also study the realization of 
certain exact complex-valued analytic $n$-forms on the standard torus through 
non-critical totally real embeddings.
\end{abstract}
\maketitle
\setcounter{thm}{0}\setcounter{equation}{0}
\section{Introduction}

In this paper we study unimodular invariants of an immersed or an embedded $n$-dimensional real manifold $M$ in complex space ${\Bbb C}^n$ with respect to 
the  holomorphic $n$-form
$\Omega=  dz_1\wedge\ldots\wedge dz_n$. We shall consider an immersion or an embedding $\varphi\colon M\to {\Bbb C}^n$ which is {\it totally real}, i.~e.
\be{omega}
\omega_\varphi =\varphi^*\Omega\neq 0.
\end{equation}
 The complex $n$-form $\omega_\varphi$ is decomposed into
\be{1.2}
\omega_\varphi=e^{i\mu}\nu,\quad i=\sqrt{-1},
\end{equation}
in which $\mu\,  (\mbox{mod}\, 2\pi)$ is a real function and $\nu$ is a real $n$-form. The decomposition~\re{1.2} is uniquely determined by  the 
requirement that $\nu$ is either a 
volume form  when $M$ is an oriented manifold, or a positive $n$-form of odd kind 
when $M$ is non-orientable. We say that the totally real immersion 
$\varphi\colon M\to {\Bbb C}^n$ is {\it non-critical\/}, if 
\be{mu}
\tau=e^{i\mu}: M\to S^1
\end{equation}
is a submersion.

There  has been extensive  investigations on real submanifolds in ${\Bbb C}^n$,
 especially since the work of E.~Bishop~\cite{bishop}. 
In association with the complex tangents, the topology of an immersed submanifold in complex space  has been studied by H.~F.~Lai~\cite{lai}, S.~M.~Webster~\cite{websterjdg}, ~\cite{webstercmh}, and most recently, by F.~Forstneri\v{c}~\cite{forstnericduke}.
In~\cite{wells}, R.~O.~Wells proved that if an orientable compact manifold $M$ admits a totally real embedding in ${\Bbb C}^n$, then its Euler number $\chi(M)$ and Pontrjagin classes must vanish. 

 From   the definition of the non-critical totally real immersion, we  have the following.  
 \bp{1.1}
 Let $M$ be a connected compact smooth manifold of dimension $n$. If $M$ admits a non-critical totally real immersion in ${\Bbb C}^n,$ then  its fundamental group $\pi_1(M)$ is infinite and the Euler number $\chi(M)$ vanishes.
\ep
As a consequence, we obtain
\bc{1.2}
The sphere $S^n$ has a non-critical totally real  immersion in ${\Bbb C}^n$  if and only if $n=1.$ If a
 compact surface $M$  admits a non-critical totally real immersion in ${\Bbb C}^2,$ then $M$ is either a torus or a Klein bottle. 
\ec
We mention that the totally real embedding of the Klein bottle in ${\Bbb C}^2$ 
given by W.~Rudin~\cite{rudin} is indeed non-critical.
 In $[21, \mbox{p.}\  26]$, A.~Weinstein constructed a Lagrangian 
(whence totally real)  
immersion of  $S^n$ into ${\Bbb C}^n$ for all $n$. An explicit totally real embedding of $S^3$ in ${\Bbb C}^3$ was given by P.~Ahern and W.~Rudin~\cite{ahern}. In~\cite{forstnericduke}, Forstneri\v{c} proved that all orientable compact surfaces
admit  totally real immersions in ${\Bbb C}^2$, and that a non-orientable surface has a totally real immersion in ${\Bbb C}^2$ if and only if its genus is even.
 
Let  $M_1, M_2$ be two totally real and real
analytic   compact submanifolds  of ${\Bbb C}^n$ with 
dimension $n$. From the 
Weierstrass approximation theorem, one knows that 
$M_1$ is equivalent to $M_2$ through biholomorphic mappings  defined near $M_1$ whenever $M_1$ and $M_2$ are diffeomorphic by smooth mappings.  In this paper, we shall study when  $M_1$  is  equivalent to $M_2$ by a unimodular transformation, i.~e.~by a  biholomorphic transformation $\Phi$  defined near $M_1$, which  satisfies the  relation $\Phi^*\Omega=\Omega$. In~\cite{gongpams}, we proved that locally,  a  non-critical  totally real
$n$-dimensional analytic submanifold $M$ in ${\Bbb C}^n$  has only one  unimodular invariant when
$n\geq 2$. Globally, we shall see that one of  unimodular invariants of $M$ is  the total volume  
$$
\rho_0=\int_M\nu.
$$
Namely,  two real-valued analytic volume  elements $\nu_1, \nu_2$ on a compact manifold $M$ are equivalent by an analytic diffeomorphism of $M$ if and only if they have the same volume,  a result due to J.~K.~Moser~\cite{mosertams}.

The main results of this paper concern the unimodular invariants of tori in ${\Bbb C}^n$  which are    perturbations of the standard torus 
$$
T^n\colon |z_j| =1,\quad  1\leq j\leq n.  
$$
Let $A_r\subset{\Bbb C}^n$ be  the annulus defined by 
$$
e^{-r}< |z_j|< e^r, \quad 1\leq j\leq n.
$$
For a  holomorphic  mapping $f=(f_1,\ldots,f_n)\colon A_r\to {\Bbb C}^n$, we denote 
$$
\|f\|_r=\sup\{ |f_j(z)|;z\in A_r,  1\leq j\leq n\}.
$$ 
We have the following result.
\bt{m}
 Let $M$ be an embedding of $T^n$ in ${\Bbb C}^n$ given by a  mapping $\varphi$
which is holomorphic in $ A_{r_0}.$  Then there exists a positive number $\epsilon_0,$ which is independent of $r_0,$ such that for $n\geq 2,\  0<r_0<1,$ and
\be{f-id}
\|\varphi-\hbox{\em Id}\|_{r_0}\leq \epsilon_0 r_0^4, 
\end{equation}
$M$ is unimodularly equivalent to an embedding of $T^n$ defined by a mapping  
$$
\psi\colon (z_1,z^\prime)\to (\zeta^{-1}
g(\zeta z_1
), z^\prime) 
$$
with $z^\prime=(z_2,\ldots,z_n)$ and $\zeta=z_2\ldots z_n, $
where $g $ is  given by
$$
\frac{d}{d\theta_1}g(e^{i\theta_1})=\rho_0e^{i(\theta_1+k(\theta_1))}
$$
for a  $2\pi$-periodic real function $k$ satisfying 
$
 \int_0^{2\pi}k( \theta_1)\, d\theta_1=0. 
$
Furthermore$,$ $k( \theta_1)$ is uniquely determined by $M$ up to a  translation $ \theta_1\to  \theta_1+\pi.$
 \et

We now consider  an  immersion problem as follows. Given a fibering $\tau$ in the form \re{mu} and a real $n$-form $\nu$ on $M$, we ask whether  there exists  a non-critical totally real  immersion $\varphi\colon M\to {\Bbb C}^n$ such that $ \varphi^*\Omega=\tau\nu$. A necessary condition is that for the complex-valued $n$-form $\omega=\tau\nu$, 
\be{mean}
\int_{M}\omega=0.
\end{equation}
For $M=T^n$, we have the following.
\bt{f}
 Let $\omega=(1+a(z))\Omega$ be a complex-valued $n$-form on $T^n.$ Assume that $a(z)$ is holomorphic in $A_{r_0}$ and $\omega$ satisfies $\re{mean}.$  Then there exists $\epsilon >0$ such that for $0<r_0<1$ and $\|a\|_{r_0}\leq \epsilon r_0,$   there is a totally real and non-critical analytic  embedding $\varphi\colon T^n\to {\Bbb C}^n$ with
$\omega=\omega_\varphi.$  
\et

We organize the paper as follows. In section 2, we shall discuss fiberings 
of a compact manifold $M$ over the circle. In particular, we shall use the triviality of fibering $T^2$ over $S^1$ to show that when $M$ is the $2$-dimensional torus,  the analytic function $\mu$ in \re{mu} has no analytic invariant with respect the whole group of analytic diffeomorphisms of $T^2$, except for the number of connected components of the fibers of $\tau$. We then normalize the function $\mu$ on $T^n$ by volume-preserving analytic transformations. This normalization, proved later in section 4, is  essential  for the proof of \rt{m}. We shall finish section 2 with a
 regular homotopy classification of non-critical immersions of $S^1$ into the complex plane. In section 3, we shall  first discuss the invariants of analytic volume forms, and then give a proof for 
 \rt{m}.  Section 5 will be devoted to the proof of \rt{f}.

\setcounter{thm}{0}\setcounter{equation}{0}
\section{Invariants of non-critical immersions}
In~\cite{tischler}, D.~Tischler proved that a fibering $M$ over $S^1$ exists whenever $M$ has a smooth foliation given by a closed $1$-form. The obstruction theory for the fiberings of  $M$ over 
$S^1$ were also studied by F.~T.~Farrell~\cite{farrell} and others (see  the references in~\cite{farrell}).
In this section, we shall first discuss a fibering $M$ over $S^1$. In particular, We shall see that for 
an analytic fibering $T^2$ over $ S^1$, the only invariant is  the number of connected components of its fiber. Then we shall 
find a complete set of invariants for perturbations of the trivial fibering $T^n$ over $ S^1$
with respect to  volume-preserving bundle maps. 
 Finally, we shall   discuss the regular homotopy classification of non-critical immersions of $S^1$
in ${\Bbb C}$,
 which is based on a classical result of H.~Whitney and W.~C.~Graustein~\cite{whitney}.
\subsection{The fibering $M$ over $S^1$}
We consider a submersion  \re{mu}.  Let $\pi\colon {\Bbb R}\to S^1$ be the covering mapping with $\pi(\theta_1)=e^{i\theta_1}$.   With
$\pi_1(S^1)={\Bbb Z}$,  we put 
$$
\tau_*\pi_1(M)={\Bbb Z}\cdot d
$$ 
for some non-negative integer $d$. It is clear that $d\neq 0$; otherwise, $\tau$  has a lifting $\widetilde\tau\colon M\to R$, which contradicts that $M$ is compact and $\tau$ is a submersion. We now choose a Riemannian metric on $M$.  Let $\nabla\mu$ be the gradient of $\mu$ 
with respect to the Riemannian metric.  Since $d\mu\neq 0$ on $M$, then the vector field $\nabla\mu$   vanishes nowhere on $ M$. This  implies that $\chi(M)=0$, and hence  \rp{1.1} is proved. 

It is a classical result of C.~Ehresmann~\cite{ehresmann} that if $M, N$ are connected compact manifolds, and $p\colon M\to N$ is a submersion,  then $M$ is a fiber bundle over $N$.
Arising from the submersion \re{mu},  the fiber bundle structure on $M$  can be described 
as follows.   Let   $\varphi_t$ be the flow of $\nabla\mu$. On $M\times{\Bbb R}$,
consider a well-defined smooth function
$$
F(x,t)=\mu\circ\varphi_t(x)-\mu(x).
$$
One notes that
   $$
\frac{d}{dt}F(x,t)=\mathopen{<}\nabla\mu(\varphi_t(x), \nabla\mu(\varphi_t(x)
\mathclose{>}
$$  
has a positive lower bound.
Hence, there is a smooth function $t\colon M\times{\Bbb R}\to{\Bbb R}$ such that
$$
F(x, t(x,\alpha))=\alpha.
$$
Now, one defines a smooth family of diffeomorphisms $\psi_\alpha$ of $M$ 
by 
$$
\psi_\alpha
(x)=\varphi_{t(x,\alpha)}(x).
$$
 In fact, one readily sees that $\alpha\to\psi_\alpha$ is a homomorphism,  i.e. an ${\Bbb R}$-action on $M$ generated by $\tau$. Now for two points $p, q$ on $S^1$ with $q=e^{i\alpha}p$
and $0<\alpha\leq 2\pi$,   we have the Poincar\'e section mapping
$$
S_{p,q}=\psi_\alpha|_{\tau^{-1}(p)}\colon 
 \tau^{-1}(p)  \to\tau^{-1}(q).
$$
For $F=\tau^{-1}(1)$, we define  two 
local trivializations 
$$
\phi_j\colon F\times (S^1\setminus {\{(-1)^j\}})\to\tau^{-1}(S^1\setminus 
{\{(-1)^j\}}),
\quad j=0,1
$$
with
 $$
\phi_0(x,p)=S_{1,p}(x), \quad \phi_1(x,p)=S_{p,-1}\circ S_{1,-1}(x).
$$
Therefore,  $\tau\colon M\to S^1$ is a fiber bundle with fiber $F$. Obviously, $\psi_\alpha$ is a bundle map between two fiber bundles $\tau\colon M\to S^1$  and $e^{i\alpha}\tau\colon M\to S^1$.

By the homotopy sequence
of a fibering (see~\cite{hu}, p.~152), we know that $d$ is exactly the number of connected 
components of
fiber $F$. Now the fibering $\tau$ has a factorization
$
\tau=(\tau^{1/d})^d,
$
where the $d$-th root $\tau^{1/d}$  is a lifting of $\tau$ for the standard $d$ to $1$ covering from $S^1$ to itself defined by $z\to z^d$; and $\tau^{1/d}\colon M\to S^1$ is a fiber bundle with connected fibers. It is clear that if  $\tau^\prime, \tau^{\prime\prime}$ are two $d$-th
roots  of $\tau$, then $\tau^{\prime\prime} =\lambda\tau^{\prime} $ with $\lambda^d=1$. Using the ${\Bbb R}$-action  on $M$ generated by the fibering $\tau^\prime $, one can find   a bundle map between the fiberings $\tau^{\prime\prime} $ and $\tau^\prime$.   Furthermore, let
$\tau_j\colon M_j\to S^1$ $(j=1, 2)$ be two fiberings. Assume that fibers of 
both $M_1$ and $M_2$  have $d$ connected components, and that $\tau_1^{1/d}$ and $\tau_2^{1/d}$ are the corresponding $d$-th roots of $\tau_1$ and $\tau_2$ respectively. Then under
bundle maps,   $\tau_1$ is equivalent to $\tau_2$ if  and only if  $\tau_1^{1/d}$ is equivalent to $\tau_2^{1/d}$.  This implies that, for the purpose of the bundle classification, we may only consider a fibering $M\to S^1$ with connected fibers.

 We  shall introduce some notations. Define a universal covering 
$R^n\to T^n$ by
\be{pi}
\pi\colon(\theta_1,\ldots,\theta_n)\to (e^{i\theta_1},\ldots, e^{i\theta_n}).
 \end{equation}
 For a mapping $\phi\colon T^n\to T^m$,   a  lifting $\widetilde \phi$ of $\phi$  is a mapping from $ {\Bbb R}^n $ to ${\Bbb R}^m$  such that $\pi\circ\widetilde \phi=\phi\circ\pi$.  Obviously, 
$\widetilde\phi$  can be written as
\be{lf}
\theta_k^\prime=\sum_{l=1}^n d_{k,l}\theta_l +f_k(\theta),\quad d_{k,l}\in {\Bbb Z},\quad k=1,\ldots, n
\end{equation}
 for some real functions  $f_k$ which are $2\pi$-periodic  in each
 variable $\theta_j$. 
 Also, a mapping $\widetilde\phi\colon R^n\to R^m$ in the form \re{lf}   generates a unique mapping $\phi\colon T^n\to T^m$ such that  $\widetilde \phi$ is a lifting of $\phi$. It is clear that if $\phi$ is a diffeomorphism,  so is   $\widetilde \phi$.
Conversely,
the degree formula (see~\cite{guillemin}, p.~189) implies that if $\widetilde \phi$ is a diffeomorphism, then $ \phi$ is a local
diffeomorphism of $T^n$ with degree 
$$
d=\det(d_{k,l}).
$$
 In particular, $\phi$ is a diffeomorphism if and only if $d=\pm 1$. Finally, 
 we shall take  $ d\theta =d\theta_1\wedge\ldots\wedge d\theta_n$ as the standard volume form on $T^n$.
A diffeomorphism $\phi$ of $T^n$ is said to be {\it volume-preserving}  
whenever $\phi$ preserves $ d\theta $.

Return to a fibering \re{mu}  with $M=T^n$. Assume that the fibering has connected fibers. Then
$$
\mu(\theta)=d_1\theta_1+\ldots+d_n\theta_n+h(\theta),
$$
where $d_1,\ldots, d_n$ are relatively prime integers, and $h$ is  a $2\pi$-periodic function. By a volume-preserving bundle map, we may assume that
\be{mu+}
\mu(\theta)=\theta_1+h(\theta).
\end{equation}
 We now restrict ourselves to the case $n=2$. Then  the fibering $\tau\colon T^2\to S^1$ is a circle bundle on $S^1$.  From the bundle classification theorem (see~\cite{steenrod}, p.~97)  and  the isotropy classification of $\mbox{Diff}(S^1)$  (see~\cite{hirsch}, p.~186), one knows that there are only two inequivalent classes of circle bundles over $S^1$. Since  the total space $T^2$ is orientable, then  $\tau\colon T^2\to S^1$ is a trivial circle bundle. This implies that
there is a smooth  bundle map 
\be{trivial}
\Phi=(\tau_1, \tau_2)\colon T^2\to S^1\times S^1
\end{equation}
with $\tau=\tau_1$. 

The following result is equivalent to the fact that an analytic circle bundle $T^2\to S^1$ is analytically trivial. For the lack of references, we  shall give a proof by using the smooth trivialization $\Phi$.
\bp{d2}
 Let $\mu$ be given by $\re{mu+}$ in which  $h(\theta)$ is an analytic $2\pi$-periodic function$.$ Assume that $d  \mu$ vanishes nowhere on $T^n.$ If $n=2,$ then there is an analytic diffeomorphism $ \Psi $ of $T^2$ such that for a lifting $\widetilde \Psi,$
$ \mu\circ\widetilde\Psi(\theta )=\theta_1.$
\ep
\pf  
We shall modify $\tau_2$  in \re{trivial}  to get an analytic bundle map $\Psi$. To this end, we  put the lifting $\widetilde\Phi$ of $\Phi$ in the form \re{lf} with $n=2$. 
  By truncating the Fourier series of $f_2$, we can find a sequence of analytic functions $g_k(\theta)$ such 
that for $ \ k\to\infty$, 
 \be{ap}
 \|g_k-f_2\|,\quad \| \partial_{\theta_j}(g_k-f_2)\|\to 0,
 \end{equation}
 where $\|\cdot\|$ stands for the $L^\infty$-norm on ${\Bbb R}^n$. We now consider the mapping $\Phi_k\colon T^2\to T^2$ defined by
$$
 (e^{i\theta_1},e^{i\theta_2})\to (\tau_1(\theta),
e^{i(d_{2,1}\theta_1+d_{2,2}\theta_2
+g_k(\theta))}).
$$
We shall denote by $D\Phi$ the Jacobian matrix of $\Phi$. Obviously, $D\Phi_k\to D\Phi$ uniformly on $T^2$ as $k\to\infty$. Hence,  $
\Phi_k$ is a local diffeomorphism for  large $k$. Since $\Phi$ and $
\Phi_k$ are of the same degree, then $\Phi_k$ is a diffeomorphism for large $k$.
 Furthermore,  they also have the same first component.
Thus, for $k$ large, $\Phi_k$ is an analytic trivialization for the fibering $\tau\colon T^2\to S^1$, which gives us 
$$\mu\circ\Phi_k^{-1}(\theta)\equiv \theta_1\ (\mbox{mod}\, 2\pi).
$$ 
 Take $\Psi=\Phi_k^{-1}$. Then, for a suitable lifting $\widetilde \Psi$, one has $\mu\circ\widetilde\Psi(\theta)=\theta_1$. This completes the proof of \rp{d2}.  $\square$

  The annulus $A_r$ has
 a universal covering 
$$
S_r\subset {\Bbb C}^n\colon |\mbox{Im}\, \theta_j|< r, \quad 1\leq j\leq n,
$$
for which  the covering mapping is given by \re{pi}. 
We have the following.
\bt{h}
 Let $h$ be a  $2\pi$-periodic  holomorphic function in $S_r,$  and let $\mu$ be defined by $\re{mu+}.$  Assume that $h(\theta )$ is real for $\theta\in
{\Bbb R}^n.$ Then there is a  constant $\epsilon>0,$ which is independent of $r,$  such that for 
$0<r<1$ and 
\be{smallh}
\|h\|_r=\sup_{\theta\in S_r}\{|h(\theta)|\} \leq \epsilon r^3,
\end{equation}
  there exist a unique $2\pi$-periodic function $k(\theta_1)$ and a volume-preserving analytic transformation $ \Phi$ of $T^n$ such that for a lifting $ \widetilde  \Phi,$
 \be{k}
  \mu\circ\widetilde\Phi(\theta)=\theta_1+k(\theta_1)
\end{equation}
with $\int_0^{2\pi} k(\theta_1)d\theta_1=0.$
\et
\pf The existence of $ \Phi$ will be given in section 4 by a KAM argument. 
Here, we only verify the uniqueness of the invariant function $k$.    Let $ \psi $  be a volume-preserving  analytic transformation of $T^n$, and  $\theta^*=\widetilde\psi(\theta)$ a lifting 
satisfying
\be{nk}
 \theta_1^*+\hat k(\theta_1^*)=\theta_1 +  k(\theta_1 )
\end{equation}
for a $2\pi$-periodic function $\hat k$ with 
$\int_0^{2\pi}\hat k(\theta_1)\, d\theta_1=0$. We need to show that $\hat k=k$.
One first notice that 
$1+k^\prime$ and $1+\hat k^\prime$ are positive on $R^1$, since $d\mu\neq 0$ on $T^n$. Hence, \re{nk} implies that   the first component of $ \widetilde \psi$  can be written as
$$
\theta_1^*=\theta_1+f_1(\theta_1)
$$
for some $2\pi$-periodic function $f_1$. The rest of components of $\widetilde  \psi$ can be put into
$$
\theta_\alpha^*=\sum_{\beta=1}^n d_{\alpha,\beta}\theta_\beta+f_\alpha(\theta),\quad \ 2\leq\alpha\leq n,
$$ 
in which $ d_{\alpha,\beta}$ are integers, and  $f_\alpha$ are $2\pi$-periodic functions. 

Denote
$$
\theta_\alpha^\prime=\sum_{\beta\geq 2} d_{\alpha,\beta}\theta_\beta,\quad \alpha\geq 2.
$$
Also,   define $d^\prime$ by
$$
d^\prime h(\theta )=\sum_{\alpha =2}^n  \partial_{\theta_\alpha}h(\theta) d\theta_\alpha.
$$
Since $\widetilde\psi$ is volume-preserving, then 
\be{det1}
 1\equiv(1+f_1^\prime(\theta_1)) \det\frac{\partial(\theta_2^\prime+ f_2(\theta), \ldots \theta_n^\prime+ f_n(\theta))}{\partial(\theta_2,\ldots,\theta_n)}.
\end{equation}
Notice that 
$$
d^\prime(\theta_2^\prime+ f_2(\theta))\wedge \ldots \wedge d^\prime(\theta_n^\prime+ f_n(\theta))
=d^\prime\theta_2^\prime\wedge\ldots\wedge d^\prime\theta_n^\prime+\ldots,
$$
where the term omitted is an exact $(n-1)$-form  in the variables $\theta_2,\ldots,\theta_n$. Thus
$$
 \int_{0\leq\theta_2,\ldots,\theta_n\leq 2\pi} d^\prime(\theta_2^\prime+ f_2(\theta))\wedge \ldots \wedge d^\prime(\theta_n^\prime+ f_n(\theta))
 $$
 is the total volume of $d^\prime\theta_2^\prime\wedge
\ldots\wedge d^\prime\theta_n^\prime$, which is obviously independent of $\theta_1$.
 Averaging the right side of \re{det1} over $0\leq \theta_\alpha\leq 2\pi
$ for $2\leq\alpha\leq n$, we see that $1+f_1^\prime(\theta_1)$ is  constant. Hence,  $f_1\equiv c$. Returning to \re{nk}, we get
$$
 k(\theta_1)= c+\hat k(\theta_1+c).
$$
By the assumption, we know that  the average values of $k$ and $\hat k$ for $0\leq\theta_1\leq 2\pi$ vanish. Averaging the above over $0\leq\theta_1\leq 2\pi$, we finally get $c=0$. Therefore, $\hat k=k$. This proves the uniqueness of the  function $k$. $\square$

\subsection{The case $n=1$}
We consider an immersion $f\colon S^1\to {\Bbb C}$. Put
$$
 \omega_f=\rho(\theta_1) e^{i(d\theta_1+h(\theta_1))}d\theta_1,\quad d\in {\Bbb Z},
$$
where 
 $h$ and $\rho$ are $2\pi$-periodic functions with $\rho>0$. Obviously,   $d=d_f$ is the Hopf degree of the Gauss map of $f$.  
One can see that the immersion $f\colon 
S^1\to {\Bbb C}$ is non-critical if and only if, locally, 
 the immersion $f$ is
strictly convex in ${\Bbb C}={\Bbb R}^2$.

 We have the following result.  
\begin{thm}[Whitney-Graustein,~\cite{whitney}]\label{whitney} 
 Let $f_0$ and $f_1$ be two immersions of $S^1$ in ${\Bbb C}.$ Then $f_1$ is regularly homotopic to $f_0$ if and only if  $d_{f_0}=d_{f_1}.$ Furthermore$,$ if both immersions $f_0$ and $f_1$ are non-critical$,$ then $f_0$ and $f_1$ are also regularly homotopic through a family of non-critical immersions$.$ 
\end{thm}
\pf The first part of the theorem is proved in~\cite{whitney}. We now modify Whitney's proof to show the second part of the theorem.    Put
$$
\omega_{f_j} =\rho_j(\theta_1) e^{i(d\theta_1+h_j(\theta_1))},\quad d_{f_0}=d=d_{f_1}.
$$
 We have $d\neq 0$ and 
$$
1+\frac{1}{d}h_j^\prime(\theta_1)>0.
$$
We now consider an orientation-preserving transformation $\psi_j$ of $S^1$ given by 
$$
\psi_j\colon e^{i\theta_1}\to e^{i(\theta_1+h_j(\theta_1)/d)},\quad j=0,1.
$$
 Obviously, $ f_j\circ\psi_j^{-1}$ is regularly homotopic to $ f_j$. 
Hence, we may still denote $ f_j\circ\psi_j^{-1}$ by $f_j$. Setting
$\widetilde f_j(\theta_1)=f_j(e^{i\theta_1})$, we then have
$$
  \widetilde  f_j^\prime ( \theta_1)=  \rho_j( \theta_1)e^{id \theta_1},\quad   \rho_j( \theta_1)>0,\quad j=0,1.
$$
 Put
$$
g_t( \theta_1)
=((1-t)   \rho_0( \theta_1)+t   \rho_1( \theta_1))e^{id \theta_1},\quad 0\leq t\leq 1.
$$
Obviously, the average value of  $g_t( \theta_1)$ on $S^1$ vanishes for all $t$.
We may assume that the average value of $\widetilde f_j $ vanishes for $j=0,1$.  Let $\widetilde f_t$ be the  unique $2\pi$-periodic function  in $\theta_1$ such that its average value  
for $0\leq \theta_1\leq 2\pi$ is zero, and
$\widetilde f_t^\prime=g_t$. It is clear that   $\widetilde f_t$ is a 
non-critical
immersion for all $t$. Now $f_t(e^{i\theta_1})=\widetilde f_t(\theta_1)$ gives us a regular homotopy $f_t$ of non-critical immersions from $f_0$ to $f_1$. The proof of  Theorem~\ref{whitney} is complete. $\square$

For later use,  we remark that a non-critical
immersion $f\colon S^1\to {\Bbb C}$ is an embedding if and only if $d_f=\pm 1$.  Here, we need a result of   Whitney~\cite{whitney}, which says that 
$$
I_f\equiv d_f- \mbox{sign}\, d_f
$$ 
is the algebraic sum of double points when $f$ is an immersion with only double points in normal crossings. It is clear that the Whitney self-intersection number $I_f$ must vanish if $f$ is an embedding. Conversely, if $f$ is non-critical
and $d=\pm 1$, one can see that the Gauss map 
$$
z\to -i\frac{f^\prime(z)}{|f^\prime(z)|}
$$
is a diffeomorphism of $S^1$.  Therefore, $f$ is an embedding. 

\section{The normal form of totally real tori}
\setcounter{thm}{0}\setcounter{equation}{0}

In this section, we shall first prove a result of Moser about normalizing analytic volume elements on a compact manifold $M$. We shall also give some relevant  estimates when $M=T^n$.
Finally, 
we shall  give a proof for \rt{m} by using \rt{h}. 

\subsection{A theorem of Moser}
Consider a non-critical totally real immersion $\varphi\colon M\to{\Bbb C}^n$. The pull-back $\omega_\varphi $ is a complex-valued $n$-form on $M$. In local coordinates $x=(x_1,\ldots,
x_n)$, we put
$$
\omega_\varphi(x)=\rho(x)e^{i\mu(x)}dx_1\wedge\ldots\wedge dx_n,
$$
where $\rho(x)$ is a positive function and $\mu(x)\, (\mbox{mod}\, 2\pi)$ is a real function. Define
$$
\nu(x)=\rho(x)dx_1\wedge\ldots\wedge dx_n.
$$
When $M$ is orientable, we require that $(x_1, \ldots,x_n)$ are  the coordinates which agree with a fixed orientation. Then $\nu$ is a volume form on $M$.  When $M$ is non-orientable, $\nu$ is still globally defined. However, $\nu$   is not an $n$-form on $M$; instead, it is  a positive $n$-form of odd kind. In both cases, $\nu$ is called a {\it volume element} on $M$. 

We have the following result.
\begin{thm}[Moser~\cite{mosertams}]\label{moser}
 Let $\nu_0$ and $\nu_1$   be two  analytic volume elements
on a compact manifold $M$ with $\int_M \nu_0=\int_M \nu_1.$
Then there is an analytic diffeomorphism $\phi$ of $M$ such that $\phi^*\nu_1=\nu_0.$
\et
\pf The proof given in~\cite{mosertams} is only for smooth volume elements. However, it works equally well in the analytic case.  Let us choose an analytic Riemannian metric on $M$. Then the Hodge decomposition theorem gives us
$$
\nu_j=d\delta\beta_j+h_j,\quad j=1,2,
$$
where $h_j$ is a harmonic $n$-form. Since $\beta_j$ is of {\it top\/} degree, we can rewrite
$$
\nu_j=\Delta\beta_j+h_j.
$$
Now the regularity of the Laplace-Beltrami operator $\Delta$ implies that  $h_j$ and $\beta_j$ are analytic $n$-forms (see~\cite{derham}, p.~177). Since $\int_M \nu_1=\int_M \nu_0$, then $h_1=h_0$.  Let 
$$
\nu_t=(1-t)\nu_0 +t\nu_1. 
$$
Then we have
$$
\nu_t=d\alpha_t+h_0,
$$
where $\alpha_t=(1-t)\delta\beta_0+t\delta\beta_1$. Obviously, $\alpha_t$ is a family of  analytic $n$-forms depending analytically on the parameter $t$. 
Using $\alpha_t$, one can construct a family of analytic diffeomorphisms $\phi_t$ such that $\phi_t^*\nu_t= \nu_0 $. For the detail, we refer to  \cite{mosertams}. $\square$

The above proof  does not provide us any estimate for the mapping $\phi$. 
For the proof of \rt{m}, we shall give some estimates of $\phi$ for the case
 $M=T^n$. 

 Let $h(\theta)$ be a $2\pi$-periodic holomorphic function defined in $S_r$.  We shall  introduce a useful decomposition  
$$
h(\theta)=\sum_{j=0}^n (L_jh)(\theta),
$$
in which $L_0h$ is a constant, and $L_jh$ depends only on $\theta_1,\ldots,\theta_j$. To ensure the uniqueness of the decomposition, we require that  
\be{ljh0}
[ L_jh ]_j =0,\quad 1\leq j\leq n,
\end{equation}
in which and also in the later discussion, we use the following notations
 $$
[f]_j=\frac{1}{2\pi}\int_0^{2\pi}f(\theta)\, d\theta_j,\quad
[f]=\frac{1}{(2\pi)^n}\int_0^{2\pi}\ldots\int_0^{2\pi} f(\theta)\, d\theta.
$$
Now, the   condition \re{ljh0} implies that
$$
L_0h(\theta)+\ldots+L_jh(\theta)=[\ldots[h]_{j+1}\ldots]_n,\quad 0\leq j<n.$$
In particular, $L_0h$ is the constant term of the Fourier series of $h$. Obviously, we have
$$
\| L_0h+\ldots+L_jh\|_ r\leq \|h\|_ r.
$$
Therefore, we get
\be{ljh}
\|L_jh\|_ r\leq 2\|h\|_ r, \quad 0\leq j\leq n.
\end{equation}
 
Denote
$$
D_jh(\theta)= \partial _{ \theta_j}h(\theta).
$$
We also define $D_j^{-1}h(\theta)$ to be the unique anti-derivative of $h $ with respect to the variable $\theta_j$ which satisfies the normalizing condition
\be{i1}
  [D_j^{-1}h ]_j=0.
\end{equation}
In addition, if $h(\theta)$ satisfies  
\be{i2}
 [h ]_j=0,
\end{equation}
then $D_j^{-1}h$ is also $2\pi$-periodic, and
$$
D_jD_j^{-1}h=D_j^{-1}D_jh=h.
$$

We need the following.
\bl{i}
 Let $h $ be a  $2\pi$-periodic holomorphic function  in $S_r$ satisfying 
$\re{i2}.$  Then
\be{i3}
\| D_j^{-1}h\|_ r\leq 2\pi\|h\|_ r.
\end{equation}
\el
\pf  We fix $\xi \in S_r$ and put $\xi_j=t+is$. We also let $\theta_k=\xi_k$
for $k\neq j$.   Since $D_j^{-1}h$ is $2\pi$-periodic in $\theta_j$, then \re{i1} implies that 
$$
 \int_{\xi_j-\pi}^{\xi_j+\pi}D_j^{-1}h(\theta)d\theta_j=0.
$$
Hence, there is $t_0\in(t-\pi,t+\pi)$ such that $\mbox{Re}\, D_j^{-1}h(\theta)=0$ for $\theta_j=t_0+is$. 
We now have
$$
\mbox{Re}\, D_j^{-1}h(\xi)=\int_{t_0+is}^{t+is}\mbox{Re}\, h(\theta)d\theta_j.
$$
Thus,
we get
$$
|\mbox{Re}\, D_j^{-1}h(\xi)|\leq \pi\|h\|_ r.
$$
One can also obtain that $|\mbox{Im}\, D_j^{-1}h(\xi)|\leq \pi\|h\|_ r.$
 Therefore, \re{i3}  is verified. $\square$

Consider a transformation $\phi $  given by
\be{id}
 \theta_j^\prime=\theta_j+f_j(\theta),\quad j=1,\ldots n,
\end{equation}
where each $f_j$ is $2\pi$-periodic on $S_r$.   We first prove the following result.
\bl{inverse}
 Let $\phi\colon S_r\to {\Bbb C}^n$ be a mapping  with the form $\re{id}.$ Assume that
 \be{nf}
\|f\|_r\leq \frac{r} {4n}.
\end{equation}
Then for $0<r<1,$  one has
  \be{domain}
 \phi\colon S_{r/2}\to S_{r},\quad 
\phi^{-1}\colon S_{ r/4}\to S_{ r/2}.
 \end{equation}
\el
\pf Obviously, \re{nf} implies that $\phi\colon S_{r/2}\to S_{r}$. For the inverse mapping, we
fix $\theta\in  S_{r/4}$ and consider a mapping 
$$
T\colon\theta^\prime\to\theta-f(\theta^\prime).
$$
From \re{nf}, it is clear that $T\colon S_{r/2}\to S_{r/2}$. By Cauchy inequalities, we
have
$$
\|  D_k f_j\|_{r/2}\leq  \frac{1}{2n}.
$$
Hence
$$
\|DT\|_{r/2}\leq 1/2,
$$
 where and also in the sequel,  the operator norm $\|Df\|_r$ of a holomorphic mapping 
$$
f=(f_1,\ldots,f_m)\colon S_r\to {\Bbb C}^m
$$ 
is defined by
$$
\|Df\|_{r}=\sup_{1\leq j\leq m}\left\{\sum_{k}|  D_kf_j(\theta)|;
\theta\in S_{r}\right\}.
$$
Therefore, we know that $T\colon S_{r/2}\to S_{r/2}$ is a contraction mapping. Now the unique fixed point $\theta^\prime\in S_{r/2}$
of $T$ is precisely the inverse image $\phi^{-1}(\theta)$. The proof of \rl{inverse} is complete. $\square$

We  have  the following estimates.
\bp{moser+} 
  Let $\omega=(1+b(\theta)) d\theta $ be a real $n$-form on $R^n.$ Assume that $b(\theta)$ is a
$2\pi$-periodic holomorphic function on $S_r$ with
\be{na}
\|b\|_r\leq  \frac{r }{ 32n\pi },\quad 0<r<1.
\end{equation}
Then there is a holomorphic mapping $\phi $ in the form $\re{id}$ 
such that
$\phi({\Bbb R}^n)={\Bbb R}^n$ and $   \omega=(1+[b])\phi^*d\theta .$ Furthermore$,$ $\phi$ satisfies $\re{domain},$ and
\be{real}
\|f\|_{r}\leq   8\pi \|b\|_r.
\end{equation}
\ep
\pf 
  We  shall seek a transformation $\phi$ in the form \re{id} such that $f_j(\theta)$ depends only on $\theta_1,\ldots, \theta_j$. Then 
$\omega=(1+[b])\phi^*d\theta$ can be rewritten as
\be{real+}
(1+[b])(1+ D_1f_1(\theta))\ldots (1+ D_nf_n(\theta))
=1+b(\theta).
\end{equation}
We further require that  
\be{fj}
[f_j]_j=0. 
\end{equation}
Thus, by applying $[\,\cdot\,]_k\ (k>j)$ to \re{real+}, 
we get
$$
(1+[b])(1+D_1f_1(\theta))\ldots (1+D_jf_j(\theta))
 =1+\sum_{l=0}^jL_lb(\theta)
$$
for $1\leq j\leq n$. Hence, 
$$
D_jf_j (\theta)=\frac{L_jb(\theta)}{1+\sum_{l=0}^{j-1}L_lb(\theta)}.
$$
From the normalizing condition \re{fj}, it follows that
$$
f_j(\theta)=\frac{D_j^{-1}L_jb(\theta)}{1+\sum_{l=0}^{j-1}L_lb(\theta)}.
$$
Using \re{ljh}, \re{i3} and \re{na}, we obtain
$$
\|f_j\|_r\leq \frac{4\pi\|b\|_r}{1-2n\|b\|_r}\leq 8\pi\|b\|_r.
$$
From \re{na} again, we see  that 
  \re{nf} is satisfied. Thus, \rl{inverse}  implies that \re{domain} holds. This completes the proof of \rp{moser+}. $\square$

\subsection{Proof of  \rt{m}}
   Let $M$ and $\varphi$ be as in \rt{m}. We take 
$$r=r_0/2,\quad 0<r_0<1.
$$
Assume that
\be{epsilon}
\epsilon_0<\frac{1}{8e^2n^2}.
\end{equation}
From Cauchy inequalities, we get
$$
\|D\varphi-\mbox{Id}\|_{r}
\leq n\frac{\|\varphi-\mbox{Id}\|_{r_0}}{e^{-r_0/2}-e^{-r_0}}\leq 
\frac{2ne}{r_0}\|\varphi-\mbox{Id}\|_{r_0}\leq 2en\epsilon_0 r_0^3,
$$
in which the last inequality is obtained from \re{f-id}.
 Let $\lambda_1(z), \ldots, \lambda_n(z)$ be the eigenvalues of $D\varphi(z)$. Then the above estimate gives us
$$
|\lambda_j(z)-1|\leq 2en\epsilon_0 r_0^3.
$$
Noticing that
\be{log}
|\log(1+x)|\leq 2|x|
\end{equation}
 for $|x|\leq 1/2$, we get from \re{epsilon} that
$$
|\sum_{j=1}^n\log \lambda_j(z)|\leq 4en^2\epsilon_0 r_0^3,\quad z\in A_{r}.
$$
Now, using the inequality 
\be{exp}
|e^{w}-1|\leq e|w|
\end{equation}
 for $|w|\leq 1$,  we  obtain 
$$
|\lambda_1(z)\ldots\lambda_n(z)-1|=
|e^{\log\lambda_1(z)\ldots\lambda_n(z)}-1|\leq 4e^2n^2\epsilon_0 r_0^3
$$
for $z\in A_{r}$. Put 
$$
1+a(z)\equiv\det(D\varphi(z))=\lambda_1(z)\ldots\lambda_n(z).
$$
 Then
\be{norma}
\|a\|_{r}\leq 4e^2n^2\epsilon_0 r^3.
\end{equation}

Let  $\widetilde{\omega_\varphi}=\pi^*\omega_\varphi$.  Put
\begin{align*}
 b_1(\theta)&=\{(1+a(z))(1+\overline{a(\overline z)}
)\}^{1/2}-1,\\
 h_1(\theta)&=\frac{1}{2i}\{\log(1+a(z))-\log(1+\overline{a(\overline z)})\}
\end{align*}
 with $z=(e^{i\theta_1},\ldots,e^{i\theta_n})$ and $\overline z\equiv(e^{i\overline\theta_1},\ldots,e^{i\overline\theta_n})$. Then, $b_1$ and $h_1$ are 
holomorphic functions in $S_r$.
On ${\Bbb R}^n$, we  now have a decomposition
$$
\widetilde{\omega_\varphi}=i^n \nu_1 e^{i \mu_1}
$$
 with
$$
  \nu_1(\theta)= (1+  b_1(\theta)) d\theta,\quad 
  \mu_1(\theta)=\theta_1+\ldots+\theta_n+  h_1(\theta).
$$ 
 From \re{epsilon}, \re{log} and \re{norma}, we obtain
$$
|\log(1+a(z))|\leq 8e^2n^2\epsilon_0 r_0^3,\quad |\log(1+\overline{a(\overline z)})|\leq 8e^2n^2\epsilon_0 r_0^3
$$
for $z\in A_{r}$. Hence, from \re{epsilon} and \re{exp}, it follows that 
  \be{bh}
\|  b_1\|_{r}\leq c_0\epsilon_0 r^3,\quad
  \| h_1\|_{r}\leq  c_0\epsilon_0 r^3,
\end{equation}
in which and also in the sequel, we denote by 
$c_j\ (j=0, 1, \ldots)$  a constant which depends only on $n$, but  is larger than $1$. 

 Consider the linear  transformation 
$$
\phi_0\colon \theta_1^\prime=\theta_1+\ldots+\theta_n, \qquad \theta_j^\prime=\theta_j, \quad j\geq 2.
$$
Put $ \mu_2=\mu_1\circ\phi_0^{-1}$ in the form
$$
\mu_2(\theta)=\theta_1+h_2(\theta),\quad h_2(\theta)=h_1\circ\phi_0^{-1}(\theta).
$$
Let
$$
\nu_2=\nu_1\circ\phi_0^{-1}=(1+b_2)d\theta,\quad b_2=b_1\circ\phi_0^{-1}.
$$
Notice that the domain $S_r$ is invariant by $\phi_0$. Then \re{bh} gives us
   \be{bh2}
\|  b_2\|_{r}\leq  c_0\epsilon_0 r^3,\quad
  \| h_2\|_{r}\leq  c_0\epsilon_0 r^3.
\end{equation}

 We now choose  $\epsilon_0$ so small that   \re{bh2}  implies that  \re{na}  holds for $b=b_2$.  Let $\phi_1$ be the transformation given in \rp{moser+}. Then 
$$
(\phi_1^{-1})^*\nu_2 =i^n(1+[b_3]) d\theta .
$$
Notice that  $\phi_0$ and $\phi_1$ are volume-preserving. Hence, 
$[b_3]=[a]$.  By setting
$$
\mu_2\circ\phi_1^{-1}(\theta)=\theta_1+ h_3(\theta),
$$
 we then have 
$$
 h_3(\phi_1(\theta))=h_2 (\theta)-f_1(\theta).
$$
Notice that $\phi_1^{-1}\colon S_{r/4}\to S_{r/2}.$ Thus
$$
\|h_3\|_{r/4}\leq \|h_2\|_{r/2}+ \|f_1\|_{r/2}.
$$
Together with \re{real} and \re{bh}, we get
\be{hh}
\|  h_3\|_{r/4}\leq c_1\epsilon_0 r^3.
\end{equation}
 Now for a possibly  smaller $\epsilon_0$, \re{hh} implies that  $h_3$ satisfies \re{smallh}. Thus, by applying \rt{h}, we obtain a volume-preserving transformation $\phi_2$ such that 
$$
\phi_2^*e^{i(\theta_1+h_3(\theta))} d\theta =e^{i(\theta_1+k(\theta_1))} d\theta .
$$
Put $\phi=\phi_0^{-1}\phi_1^{-1}\phi_2 \phi_0 $. Then, we obtain
\be{nomega}
\phi^*\widetilde {\omega_\varphi}=i^n(1+[a])e^{i(\theta_1+\ldots\theta_n+k(\theta_1+\ldots\theta_n))} d\theta .
\end{equation}

To construct $g$ as stated in \rt{m}, we first notice that $\phi^*\widetilde {\omega_\varphi}$ is exact, i.e.
$$
\int_{T^n}e^{i(\theta_1+\ldots\theta_n+k(\theta_1+\ldots\theta_n))} d\theta =0.
$$
Equivalently, this  means that the average value of $e^{i(\theta_1+k(\theta_1))}$ for $0\leq \theta_1\leq 2\pi$ vanishes.  Hence, there is an immersion
$$
g\colon S^1\to {\Bbb C}
$$
such that
$$
\frac{d}{d\theta_1}g(e^{i\theta_1})=e^{i(\theta_1+k(\theta_1))}.
$$
 We further notice that $\varphi\circ\phi$ is a non-critical totally real immersion. Hence, from the right side of \re{nomega}, we see that $1+k^\prime(\theta_1)\neq 0$. Therefore, $g$ is a non-critical
immersion 
with $d_{g}=1$. 
As we have seen in the end of section 2, this  implies that $g$ is an embedding.
Now, one can verify that the mapping $\psi$ given in \rt{m} defines an embedding of $T^n$ with
$
\psi^*\Omega=\phi^*\varphi^*\Omega.
$
 Let $\widetilde T^n=\psi(T^n)$. Then the analytic diffeomorphism $\Phi=\psi\circ\phi^{-1}\circ\varphi^{-1}\colon M\to 
\widetilde T^n$ satisfies
$$
\Phi^*\Omega|_{\widetilde T^n}=\Omega|_M.
$$
By complexifying $\Phi$, we get a unimodular transformation defined near $M$  which transforms $M$ into $\widetilde T^n$. 

Finally, we need to verify that up to a translation 
$\theta_1\to \theta_1+\pi$, the function $k(\theta_1)$ is uniquely determined by $M$. More precisely, we want to show that if $M$ is unimodularly equivalent to an analytic embedding of $T^n$ given by another mapping
$$
 \hat \psi\colon (z_1,z^\prime)\to (\zeta^{-1} \hat g(\zeta z_1), z^\prime)
$$
with
$$
\frac{d}{d\theta_1} \hat g(e^{i\theta_1})= \hat \rho_0e^{i \hat k(\theta_1)},\quad [\hat k]=0.
$$
Then, either $ \hat k=k$, or $ \hat k(\theta_1+\pi)=k(\theta_1).$  For the proof, we notice that
$$
\omega_{ \hat \psi}= \hat\rho_0e^{i(\theta_1+\ldots+\theta_n+ 
\hat k(\theta_1+\ldots+\theta_n))}.
$$
If $M$ is unimodularly equivalent to $ \hat \psi(T^n)$, then there is a diffeomorphism $\phi$ of $T^n$ such that  $\phi^*\omega_{ \hat \psi}=\omega_\psi$. In brief,  let us still denote the lifting of $\phi$ by $\phi$. This implies that either
\be{case1}
 \hat \rho_0\phi^* d\theta =\rho_0 d\theta 
\end{equation}
and 
\be{case1+}
 \hat K\circ\phi(\theta)=K(\theta)+2d\pi,\quad
d\in{\Bbb Z}
\end{equation}
for $K(\theta)=\theta_1+k(\theta_1)$ and $ \hat K(\theta)=\theta_1+ \hat k(\theta_1)$; or
\be{case2}
 \hat \rho_0\phi^* d\theta =-\rho_0 d\theta 
 \end{equation}
and
\be{case2+}
 \hat K\circ\phi(\theta)=K(\theta)+(2d+1)\pi.
\end{equation}
If \re{case1} holds, then we  see that  $ \hat \rho_0=\rho_0$ and $\phi$ is volume-preserving.   Let $L_a$ be the translation  
\be{la}
\theta_1\to \theta_1+a, \qquad \theta_j\to \theta_j, \quad j\geq 2.
\end{equation}
Notice that $k$ is $2\pi$-periodic. Then \re{case1+} implies that $ \hat K\circ\phi\circ L_{-2d\pi}=K$.
Since $\phi\circ L_{-2d\pi}$ is volume-preserving, then \rt{h} implies that $ \hat k=k$. 
 Next, we assume that \re{case2} holds. 
In this case, we get $ \hat \rho_0=\rho_0$ and $\phi^* d\theta =- d\theta .$ 
Here we use the assumption $n\geq 2$.  We  consider the  volume-preserving mapping
$$
\hat\phi(\theta)=\phi(\theta_1,\ldots,\theta_{n-1},-\theta_n).
$$
Notice that $ \hat K$ depends only $\theta_1$. Hence, we still have
$$
 \hat K\circ\hat\phi(\theta)=K(\theta)+(2d+1)\pi.
$$
Now \re{case2+} can be rewritten as $ \hat K\circ\hat\phi\circ L_{(2d+1)\pi}=\theta_1+k(\theta_1+\pi)$. Again, \rt{h} 
implies that $ \hat k(\theta_1)=k(\theta_1+\pi)$. This completes the proof of \rt{m}.

\setcounter{thm}{0}\setcounter{equation}{0}
\section{Normalizing functions  with volume-preserving mappings}
 
This section is devoted to the proof of \rt{h}. We shall construct volume-preserving transformations 
of $T^n$ through the time-1 mappings of divergent free vector fields, since the latter  are easier to be handled because of the linearity of the compatibility condition. We shall accomplish the proof of \rt{h} by using a KAM argument.

\subsection{Some estimates for flows}
Consider the flow $\phi_t$ determined by a system of ordinary differential equations
\be{ode}
\frac{d\theta_j}{dt}=p_j(\theta), \quad 1\leq j\leq n.
\end{equation}
Set
$$
\phi_t: \theta_j^\prime=\theta_j+f_j(\theta,t), \quad 1\leq j\leq n.
$$
For $p=(p_1, \ldots,  p_n)$, we put
$$
\|p\|_r=\sup \{|p_j(\theta)|;\theta\in S_r,1\leq j\leq  n\}.
$$

We need the following result. 
\bl{z}
  Let $ p_1,\ldots, p_n$ be  holomorphic functions on $S_{r_1}.$ 
 Assume that
\be{z1}
\|p\|_{r_1}\leq  r_1\delta.
\end{equation}
Then for $-1\leq t\leq 1$ and $0<\delta<1/2,$ one has
\begin{gather}\label{z2}
\|f(\cdot,t)\|_{(1-\delta)r_1}\leq  \|p\|_{r_1}\\
\intertext{and}
\|f_j(\cdot, t)-p_j(\cdot)t\|_{(1-2\delta) r_1}\leq\frac{n\|p\|_{r_1}^2}{ r_1\delta}.\label{z3}
\end{gather}
Furthermore$,$  if all $p_j$ are $2\pi$-periodic, then $f_j(\theta,t)$ are also $2\pi$-periodic in 
$\theta$ for  fixing $|t|\leq 1.$
\el
\pf  We first want to prove that $\phi_t(\theta)$ is holomorphic for $\theta\in S_{(1-\delta)r_1}$ and $|t|\leq 1$. By the Cauchy existence theorem, it suffices to show that for any $t$ with   $|t|\leq 1$, one has 
\be{odedomain}
\phi_t(S_{(1-\delta)r_1})\subset S_{r_1}. 
\end{equation}
Otherwise, there are $0<|t_0|\leq 1$ and $\theta_0\in S_{(1-\delta)r_1}$ such that
for some $k$ 
$$
|f_k(\theta_0,t)|<|f_k(\theta_0,t_0)|=\delta r_1,\quad |t|<|t_0|.
$$
Rewrite \re{ode} as
 \be{node}
f_j(\theta, t)=\int_0^tp_j(\phi_\tau(\theta))d\tau,\quad j=1,\ldots, n.
\end{equation}
Then \re{z1} implies that
$$
|f_k(\theta_0,t_0)|<\|p_k\|_{r_1}\leq r_1\delta,
$$
which leads to a contradiction. Hence,  we have verified that \re{odedomain} holds. Using \re{node}, 
we now obtain the inequality (\ref{z2}) for $|t|\leq 1$.
  
 By Cauchy inequalities, we have
\be{z5}
\|D_lp_k \|_{(1-\delta) r_1}\leq\frac{\|p_k\|_{r_1}}{ r_1\delta}.
\end{equation}
Fix $\theta\in S_{(1-2\delta) r_1}$ and $|t|\leq 1$. Then (\ref{z2}) 
and \re{odedomain} give us
$$
|p_k(\theta+f(\theta,t))-p_k(\theta)|\leq
\sum_{l=1}^n\|f_l(\cdot, t)\|_{(1-2\delta) r_1}\cdot\| D_lp_k(\cdot, t) \|_{(1-\delta) r_1}.
$$
Together with (\ref{z2}) and \re{z5}, we get
$$
|p_k(\theta+f(\theta,t))-p_k(\theta)|\leq\frac{n\|p\|_{r_1}^2}{ r_1\delta}.
$$
Now, from the identity
$$
f_k(\theta,t)-p_k(\theta)t=\int_0^t\{p_k(\theta+ f(\theta,\tau))-p_k(\theta)\}d\tau,
$$
we  obtain the estimate (\ref{z3}). The periodicity of $f_j$ follows directly from the uniqueness of the solutions $f_1,\ldots,f_n$ to  \re{ode} and the assumption that all $p_k$ are $2\pi$-periodic. The proof of \rl{z} is complete.
$\square$

We notice 
that if the vector field $v$ defined by \re{ode} is {\it divergent free}, i.e.
\be{div}
\sum_{j=1}^n  D_jp_j(\theta)=0,
\end{equation}
then $\phi_t$ is volume-preserving for all $t$. 

\subsection{Approximate equations}
We put
\be{p2}
p_1=\frac{1}{1+D_1L_1h}\sum_{j=2}^nL_jh.
\end{equation}
In order  that \re{div} is satisfied,  we take
\be{p3}
p_j=D_j^{-1}D_1\left(\frac{-L_jh}{1+D_1L_1h}\right),\quad j=2,\ldots,n.
\end{equation}
Let $\phi_t$ the flow defined by \re{ode}.  Set
\be{p7}
\phi=\phi_{-1}\colon \theta_j^\prime=\theta_j +f_j(\theta ,-1).
\end{equation}
With the above transformation, we now put
\be{k+}
\theta_1 +k(\theta )=\theta_1^\prime+h(\theta^\prime).
\end{equation}

Introduce the notations
$$
 B_r=\max\{ |L_0h|,\|D_1L_1h\|_r\},\quad b_r=\max \{ \| L_ 2h\|_ r,\ldots, \| L_ nh\|_ r\}.
$$
We need the following.
\bl{p}
 Let $p_j$ be defined by $\re{p2}$ and $\re{p3}.$ Assume that
\be{p4}
B_ r\leq 1/2.
\end{equation}
Then there exists a constant $c_2$ 
such that for $0<r<1$ and $0<\delta<1/2,$ one has
\be{p5}
\|p\|_{(1-\delta) r}\leq \frac{c_2b_ r}{ r\delta}.
\end{equation}
\el
\pf Using Cauchy inequalities, we get
 \begin{gather}\label{dh1}
\|D_1^2L_1h\|_{(1-\delta) r}\leq \frac{B_r}{ r\delta},\\
\|D_1L_jh\|_{(1-\delta) r}\leq \frac{b_ r}{ r\delta},\quad 2\leq j\leq n.\label{dh2}
 \end{gather}
From \re{p4} and \re{p2},  it follows that
\be{p6}
\|p_1\|_{(1-\delta) r}\leq 2nb_ r.
\end{equation}
We also have 
$$
D_1\left(\frac{L_jh}{1+D_1L_1h}\right)=\frac{1}{(1+D_1L_1)^2}\left\{(1+D_1L_1h)D_1L_jh-D_1^2L_1h \cdot L_jh\right\} .
$$
 Now, \re{i3}, (\ref{dh1}) and (\ref{dh2}) give us 
$$
\|D_j^{-1}D_1\left(\frac{L_jh}{1+D_1L_1h}\right)\|_{(1-\delta) r}\leq 8\pi\left((1+B_r)
\frac{b_ r}{ r\delta}+\frac{B_r b_ r}{r\delta}\right)\leq 
 \frac{16\pi b_ r}{ r\delta},
$$
where the last inequality comes from \re{p4}.
Together with \re{p6}, we can obtain \re{p5} by taking $c_2=\max\{2n, 16\pi\}$.
$\square$

To apply \rl{z}, we  should take 
$$ r_1=(1-\delta) r, \quad 0<r<1, \quad 0<\delta<1/4.
$$ 
We also
assume that
\be{b}
b_ r\leq \frac{ r^2\delta^2}{nc_2}.
\end{equation}
Notice that $ r_1> r/2$. Then, \re{p5} and \re{b} imply that \re{z1} holds. 
 Since $(1-2\delta) r<(1-\delta)r_1$,
 then \re{z1} and (\ref{z2}) give us
\be{it3}
\|f(\cdot, -1)\|_{(1-2\delta)r}\leq r\delta.
\end{equation}
From (\ref{z3}) and \re{p5}, we also have
\be{p9}
\|f_j(\cdot,-1)+p_j(\cdot)\|_{(1-3\delta) r}\leq \frac{nc_2^2b_ r^2}{ r_1 r^2\delta^3} 
\leq
 \frac{c_3b_ r^2}{ r^3\delta^3}.
\end{equation}

Rewrite \re{k+} as
\be{p12}
k(\theta)=h(\theta-p(\theta))-p_1(\theta)+E_1(\theta)
\end{equation}
with
$$
E_1(\theta)=f_1(\theta,-1)+p_1(\theta)+h(\theta+f(\theta,-1))-h(\theta-p(\theta)).
$$
Using the Taylor formula, we get
\be{p13}
h(\theta-p(\theta))=h(\theta)-\sum_{j=1}^n D_jh(\theta)\cdot p_j(\theta)+E_2(\theta)
\end{equation}
with
$$
E_2(\theta)= \sum_{i\leq j} p_i(\theta)p_j(\theta) D_iD_jh(\theta-tp(\theta)) 
$$
for some $t\in (0,1)$.
Notice that $L_1h$  is a function in $\theta_1 $ alone. Then we can write
\be{p14}
\sum_{j=1}^n  D_jh(\theta)p_j(\theta)
=D_1L_1h(\theta_1)p_1(\theta)+E_3(\theta)
\end{equation}
with
$$
E_3(\theta)=\sum_{j=2}^n\sum_{k=1}^jD_kL_jh(\theta)p_k(\theta).
$$
 Substituting \re{p13} and \re{p14}  into \re{p12}, we get
$$
k =h -(D_1L_1h+1)p_1+E,\quad E=E_1 +E_2 -E_3 .
$$
Using \re{p2}, we obtain
\be{p15}
k(\theta)=L_0h+L_1h(\theta_1)+E (\theta) .
\end{equation}

We now give some estimates of $E_j$.
Using Cauchy inequalities,  one gets
\be{p16}
\|D_jh\|_{(1-\delta) r}\leq B_r+(n-1)\frac{b_r}{r\delta}<1,
\end{equation}
in which \re{p4} and \re{b} are used. One also has
\be{p16+}
\|D_iD_jh\|_{(1-\delta) r}\leq \frac{1}{ r\delta}.
\end{equation}
Fix  $\theta\in S_{(1-3\delta) r}$. From \re{z1}, we see that $\theta-p(\theta)\in S_{(1-2\delta) r}$. Hence
$$
|h(\theta+f(\theta,-1))-h(\theta-p(\theta))|\leq
 \|Dh\|_{(1-2\delta) r}\cdot\|f(\cdot, -1)+p(\cdot)\|_{(1-3\theta) r}. $$ 
  Together with \re{p9} and \re{p16}, we get
$$
|h(\theta+f(\theta,-1))-h(\theta-p(\theta))|\leq \frac{c_3n b_ r^2}{ r^3\delta^3}.
$$
 Hence
\be{p17}
\|E_1\|_{(1-3\delta) r}\leq \frac{c_3(n+1) b_ r^2}{ r^3\delta^3}.
\end{equation}
Notice that $\theta-tp(\theta)\in S_{(1-2\delta) r}$ for $|t|\leq 1$. From \re{p5} and \re{p16+}, it follows that
\be{p18}
\|E_2\|_{(1-3\delta) r}\leq \frac{n^2}{ r\delta}\left(\frac{c_2 b_ r}{ r\delta}\right)^2.
\end{equation}
It is clear that for  $j\geq 2$,  $\|D_kL_jh\|_{(1-3\delta) r}\leq b_ r/(3 r\delta)$. Together with \re{p5}, we get
\be{p19}
\|E_3\|_{(1-3\delta) r}\leq \frac{n(n-1)b_ r}{ 3r\delta}\cdot\frac{c_2b_ r}{ r\delta}.
\end{equation}
 Let us put \re{p17}, \re{p18} and \re{p19} together in the form
\be{p20}
\|E_1\|_{(1-3\delta) r}+\|E_2\|_{(1-3\delta) r}+\|E_3\|_{(1-3\delta) r}\equiv E_0\leq
\frac{c_4b_ r^2}{2 r^3\delta^3}. 
\end{equation}

We now denote
 \begin{align*}
  \widetilde B_{(1-4\delta) r}&=\max\{\|L_0k\|_{(1-4\delta) r},
\|D_1L_1k\|_{(1-4\delta) r}\},
 \\
\widetilde b_{(1-4\delta) r}&=\max\{ \|L_2k\|_{(1-4\delta) r},\ldots,
\|L_nk\|_{(1-4\delta) r }\}.
 \end{align*}
 Applying $D_1L_1$ to both sides of  \re{p15},  we get  
$$
\|D_1L_1k\|_{(1-4\delta) r}\leq \|D_1L_1h\|_{(1-4\delta) r}+
\|D_1L_1 E \|_{(1-4\delta) r}.
$$
 From \re{ljh} and Cauchy inequalities, we then obtain
 \be{p21+}
\|D_1L_1k\|_{(1-4\delta) r} \leq B_ r+\frac{2}{r\delta}E_0.
\end{equation}
It is easy to see that $L_iL_j\equiv 0$ for $i\neq j$. Hence, applying $L_0$ to both sides of  \re{p15}
gives
\be{p21++}
|L_0k|\leq B_r+E_0.
\end{equation}
We also  have
$$
\|L_jk\|_{(1-3\delta) r}\leq  \|L_jE\|_{(1-3\delta) r},\quad j> 1.
$$
Hence
\be{p22}
\widetilde b_{(1-3\delta) r}\leq  2E_0.
\end{equation}

Substituting \re{p20} into \re{p21+}-\re{p22}, we obtain the following. 
\bp{it}
Assume that $0<r<1$ and $0<\delta<1/4.$ Suppose that $h$ is a $2\pi$-periodic holomorphic function  satisfying $\re{p4}$ 
and $\re{b}.$ Let $\phi$ be defined by $\re{p7},$ and let $k$ be determined  by $\re{k+} $ and $\re{p7}.$
Then $\phi$ satisfies $\re{it3},$ and 
 $k$  satisfies 
\be{it1}
 \widetilde  b_{(1-4\delta) r}\leq \frac{c_4b_ r^2}{ r^3\delta^3}
\end{equation}
and
\be{it2}
 \widetilde B_{(1-4\delta) r}\leq B_ r+\frac{c_4b_ r^2}{ r^4\delta^4}.
\end{equation}
\ep
 
\subsection{Proof of \rt{h}}

We put
$$
r_m=\frac{1}{2}\left(1+\frac{1}{m+1}\right)r,\quad m=0,1,\ldots.
$$
 Rewrite
 $$
r_{m+1}=(1-4\delta_m)r_m,\qquad \delta_m=\frac{1}{4(m+2)^2},\quad m=0,1,\ldots.
$$
 Let us first prove a numerical result.
\bl{4.2}
 Let $r_m,\delta_m$ be given as above$.$ Assume that $b_m$ and $ B_m\ (0\leq m\leq N)$ are
 non-negative numbers 
satisfying
\be{4.7}
b_{m+1}\leq
\frac{c_4b_m^2}{r_m^3\delta_m^3 }
\end{equation}
and
\be{4.7+}
B_{m+1}\leq B_m+\frac{c_4b_m^2}{r_m^4\delta_m^4}
 \end{equation}
for $0\leq m<N.$ Assume further that 
\be{4.8}
B_0\leq 1/4,\quad b_0\leq r_0^3\delta_0^3/c_6
\end{equation} 
with $c_6=\max\{nc_2,  27c_4\}.$ Then for all $m>0,$ we have
\be{4.6}
b_m\leq r_m^3\delta_m^3/c_6,\quad B_m\leq 1/2.
\end{equation}
\end{lemma}
\pf
We put 
$$
\hat b_m=r_m^3\delta_m^3/c_6,\quad m\geq 0.
$$
Assume that \re{4.6} holds for $m\leq m_0$. It is easy to see that
 $$
\hat b_{m+1}/\hat b_m=(1-4\delta_m)^3\frac{ \delta_{m+1}^3}{ \delta_m^3}\geq(1-4\delta_0)^3\frac{ \delta_{1}^3}{ \delta_0^3}
 =1/ 27.
$$
 On the other hand, one has
$$
\frac{c_4b_{m_0}}{r_{m_0}^3\delta_{m_0}^3}\leq c_4/c_6\leq 1/ 27.
$$
Hence
\be{4.9}
b_{m_0+1}\leq b_{m_0}\frac{\hat b_{m_0+1}}{\hat b_{m_0}}\leq \hat b_{m_0+1}.
\end{equation}
As for the estimate of $B_{m_0+1}$, we have
$$
B_{m_0+1}\leq B_0+\sum_{j=0}^{m_0} \frac{c_4b_j^3}{r_j^4\delta_j^4}\leq 1/4+\sum_{j=0}^{m_0}r_j^5\delta_j^5,
$$ 
in which the second inequality is obtained from the estimate of $b_m$ for $m\leq m_0+1$.
Notice that $r_0<1$. Then it is easy to see that $B_{m_0+1}\leq 1/2$. Therefore, we have proved \re{4.6} by the induction. $\square$

For the proof of \rt{h}, we shall find a sequence of volume-preserving transformations $\phi_m$, and then define
$$
\theta_1 +h^{(m+1)}(\theta )=\theta_1^\prime+h^{(m)}(\theta^\prime),\quad h^{(0)}
=h
$$
with
$$
\phi_m\colon \theta^\prime=\theta+f^{(m)}(\theta,-1).
$$
Denote
\begin{gather*}
B_m=\max\{|L_0h^{(m)}|, \|D_1L_1h^{(m)}\|_{r_m}\},\\
 b_m=\max\{ 
\|L_2h^{(m)}\|_{r_m},\ldots,\|L_nh^{(m)}\|_{r_m}\}.
\end{gather*}
The transformations $\phi_m$ will be constructed such that
\be{it3+}
\|f^{(m)}(\cdot, -1)\|_{r_{m+1}}\leq r_m\delta_m.
\end{equation}
In particular,
\be{m}
\phi_m\colon S_{r_{m+1}}\to S_{r_m}.
\end{equation}
Furthermore, it required that \re{4.6} holds.

The existence of $\phi_0$ is based on the initial condition \re{smallh}. 
In fact, we should choose $\epsilon$ such that \re{4.8} follows from 
\re{smallh}. Hence, we get \re{p4} and \re{b} for $h=h^{(0)}$. Let 
$\phi_0=\phi$ be as in \rp{it}. Then \re{it3} implies that
\re{it3+} holds for $m=0$. Inductively, we assume that $\phi_0,\ldots, \phi_{m_0}$ are constructed such that \re{it3+} and
\re{4.6} are satisfied for $m\leq m_0$. Then \rp{it} implies that there exists a transformation $\phi^{(m+1)}$ satisfies \re{it3+}. 
 The existence of $\phi_0,\ldots, \phi_{m_0+1}$ also implies that \re{it1} and \re{it2} hold for $m\leq m_0+1$.  Finally, from \rl{4.2}, we obtain \re{4.6} for $m=m_0+1$.  By the induction, we may construct a sequence of transformations $\phi_m$ satisfying \re{it3+} such that the corresponding $h^{(m)}$ satisfies \re{4.6}. 

We now want to prove that $\phi_0\circ\ldots\circ\phi_m$ converges on
$S_{r_0/2}$ as $m\to\infty$. We first notice that $S_{r_0/2}\subset S_{r_m}$ for all $m$.  
From \re{m}, it follows that
\begin{align*}
\|\phi_0\circ\ldots\circ\phi_m-\hbox{Id}\|_{r_0/2}&\leq 
 \sum_{j=0}^m \|\phi_j\circ\ldots\circ\phi_m- 
\phi_{j+1}\circ\ldots\circ\phi_m\|_{r_0/2}\\
&\leq\sum_{j=0}^m \|f^{(j)}(\phi_{j+1}\circ\ldots\circ\phi_m, -1)\|_{r_0/2}.
\end{align*}
On the other hand, from \re{it3+}, we get
$$
\|f^{(m)}(\cdot, -1)\|_{r_{m+1}}\leq \delta_m.
$$
Therefore,
 $\phi_0\circ\ldots\circ\phi_m$ converges to a holomorphic mapping
 $\Phi_\infty \colon S_{r_0/2}\to {\Bbb C}^n$.  It is clear that $\Phi_\infty-\mbox{Id}\colon{\Bbb R}^n\to {\Bbb R}^n$ is $2\pi$-periodic.  Hence, the restriction of $\Phi_\infty$ to ${\Bbb R}^n$ generates a volume-preserving diffeomorphism $\Phi$ of $T^n$ satisfying \re{k}, where
$$
k=\lim_{m\to\infty}\,\{L_0h^{(m)}+L_1h^{(m)}\}.
$$
 One can further achieve $[k]=L_0k=0$ by applying a suitable translation \re{la}. 
Therefore, the proof of \rt{h} is complete.

\setcounter{thm}{0}\setcounter{equation}{0}
\section{The Realization of complex $n$-forms}

In this section, we shall give a proof for \rt{f}.
As in section 4, we shall first 
consider an approximation to the nonlinear equation. We then use the KAM method to show the existence of a convergent solution to the original equation.
\subsection{A linearized equation} For the proof of \rt{f}, we try to find a time-1 mapping $\varphi$ of a holomorphic vector field in $A_r$
given by 
\be{v}
v(z)= \sum_{j=1}^n q_j(z) \frac{\partial}{\partial z_j}
\end{equation}
  such that $\varphi^*\Omega=\omega$, i.~e.
$$
\det D\varphi =1+a(z).
$$
Notice
the following determinant formula (see~\cite{siegelbook}, p.~142)
\be{d}
\log \det D\varphi_t =\int_0^t\sigma\circ\varphi_s(z)\, ds,\quad -1\leq t\leq 1
\end{equation}
with
$$
\sigma(z)=\sum_{j=1}^n  \partial_jq_j,\quad  \partial_j\equiv \frac{\partial}{\partial z_j}.
$$
Thus, we are led to the following functional equation
\be{f1}
\int_0^1\sigma\circ\varphi_t(z)\, dt=\log(1+a(z)). 
 \end{equation}

Instead of solving \re{f1}, we now consider its linearized equation 
\be{f2}
\sum_{j=1}^n \partial_jq_j=a(z).
\end{equation}
 To solve the linearized equation, we decompose
$$
a(z)=K_1a(z)+ \ldots+ K_na(z)+ K_{n+1}a(z)
$$
 with
$$
  K_ja(z)=\frac{1}{z_1\ldots z_{j-1}}\sum_{  i_j\neq -1} a_{-1,\ldots, -1,i_j,\ldots, i_n}z_j^{i_j}\ldots z_n^{i_n}
$$
for $1\leq j\leq n$. It is easy to see that   the condition \re{mean} holds for $\omega$,
if and only if 
\be{kn}
K_{n+1}a(z)=\frac{a_{-1,\ldots, -1}}{z_1\ldots z_n}\equiv 0.
\end{equation}
We also have
$$
\sum_{l=j}^{n+1} K_la(z)  =\frac{z_1^{-1}\ldots z_{j-1}^{-1}}{(2\pi i)^{j-1}}
\int_{|z_1|=\ldots=|z_{j-1}|=1}
a(z)\,dz_1\wedge\ldots\wedge z_{j-1}.
$$
Notice that $|z_j|<e$ for  $z\in A_r$, $0<r<1$.
Hence
\be{f2+}
\| K_ja\|_r\leq 2e^j\|a\|_r.
\end{equation}

 We now let $q_j(z)$ be the unique Laurent series, containing no terms of the form 
$z_1^{i_1}\ldots z_n^{i_n}$ with $i_j=0$, such that 
$\partial_jq_j(z)=K_ja(z).$
 Then \re{kn} implies that $q_1,\ldots, q_n$ satisfy \re{f2}. 
It is easy to see that
\be{f3}
q_j(z)= i D_j^{-1}(z_jK_ja(z)),\quad z=(|z_1|e^{i\theta_1},\ldots,|z_n|e^{i\theta_n}).
\end{equation}
 Then, it follows from \re{i3} and \re{f2+} that
\be{f3+}
\|q_j\|_r\leq 4\pi e^{j+1}\|a\|_r,\quad 0<r<1.
\end{equation}

 Let $v$ be the   vector field defined by
 $\re{v}$ and  $\varphi_t$  the flow of $v.$ Put 
$$
\varphi_t\colon z_j^\prime=z_jg_j(z,t),\quad g_j(z,0)\equiv 1.
$$
We have the following.
\bl{f1}
Let $a(z), \varphi_t, g_j(z,t)$ be as above.  Assume that  
\be{f4}
\|a\|_r\leq \frac{r\delta^2}{4e^{n+2}\pi}.
\end{equation}
 Then for $0<r<1, \  0<\delta<1/2$ and $|t|\leq 1,$ we have
\be{f5}
\|\log g_j(\cdot, t)\|_{(1-\delta)r}\leq r\delta^2.
 \end{equation}
  \el
\pf  We consider the vector field 
$$
\widetilde v(\theta)=\sum_{j=1}^np_j(\theta)\frac{\partial}{\partial \theta_j}, \quad 
 p_j(\theta)=-ie^{-i\theta_j}q_j(\pi(\theta)).
$$
Clearly, $\widetilde v$ is a holomorphic vector field in $S_r$ with  $\pi_*\widetilde v=v$.  
Let $\phi_t$ be the flow of $\widetilde v$. Then $\varphi_t=\pi\phi_t\pi^{-1}.$ Hence
\be{2flows}
g_j(z,t)=e^{if_j(\log z_1,\ldots, \log z_n,t)}.
\end{equation}
Since $0<r<1$, we have
$
\|p_j\|_r\leq e\|q_j\|_r.
$  
From \re{f3+} and \re{f4}, it is easy to see that 
\be{f4+}
\|p_j\|_r\leq r\delta^2.
\end{equation}
In particular, $p$ satisfies \re{z1} for $r_1=r$. 
 Now \re{f5} is obtained from (\ref{z2}), \re{2flows} and \re{f4+} immediately.
 This completes the proof of \rl{f1}. $\square$

Let $\varphi_t$ be the flow given in \rl{f1}. We put $\psi=\varphi_{-1}$ and denote
$$
\psi^*\omega=\hat \omega, \quad  \hat\omega(z)=(1+\hat  a(z))dz_1\wedge\ldots\wedge dz_n.
$$
Then 
\be{hata}
 1+\hat a(z)=(1+a\circ\psi)\det D\psi.
\end{equation}
We are ready to prove the following.
\bp{f2}
 Let $a(z)$  and $\varphi_t$ be given as in $\rl{f1}.$  Let $\hat  a(z)$ be  defined as above$.$ Then  
\be{varphit}
\varphi_t\colon A_{(1-2\delta)r}\to A_{(1-\delta)r}
\end{equation}
and
 \be{a}
\|  \hat a\|_{(1-2\delta)r} \leq \frac{c_7\|a\|_r^2}{r\delta}
\end{equation}
for some constant $c_7$   depending only on $n.$ 
\ep
\pf From \re{f5}, it follows that \re{varphit} holds. 
Let
$$
Q=\log(1+a\circ\psi)+\log\det D\psi.
$$
From \re{d},  we  get
\be{f7}
Q
=\log(1+a\circ\psi)-a\circ\psi+
\int_0^{-1}(a\circ\varphi_s-a\circ\psi)\, ds.
\end{equation}
By Cauchy inequalities, we have
$$
\|\partial_j a\|_{(1-\delta)r}\leq \frac{e\|a\|_r}{r\delta}.
$$
Notice that
$$
a(\varphi_t(z))-a(\psi(z))=\sum_{j=1}^n\int_{-1}^t \partial_ja(\varphi_s(z))q_j(\varphi_s(z))\, ds.
$$
Hence, \re{varphit} leads to
$$
\|a\circ\varphi_t-a\circ\psi\|_{(1-2\delta)r}
\leq \frac{ne}{r\delta}\|a\|_r\|q\|_r.
$$
From \re{f3+},  we then obtain
$$
\|\int_0^{-1}(a\circ\varphi_s-a\circ\psi)\, ds\|_{(1-2\delta)r}\leq \frac{4n\pi e^{n+2}}{r\delta}\|a\|_r^2.
$$
Notice that $\|a\|_r<1/2$. Hence, \re{log} gives us
$$
\|\log(1+a\circ\psi)-a\circ\psi\|_{(1-2\delta)r}\leq 2\|a\|_r^2.
$$
From \re{f7}, we now see that
\be{f8}
\|Q\|_{(1-2\delta)r}\leq
\frac{4n\pi e^{n+2}+1}{r\delta}\|a\|_r^2.
\end{equation}
In particular, \re{f4} and  \re{exp}
imply that
$$
\|\hat a\|_{(1-2\delta)}=\|e^Q-1\|_{(1-2\delta)}\leq  e\|Q\|_{(1-2\delta)}.
$$
Together with \re{f8}, we obtain \re{a}. 
$\square$

\subsection{Proof of \rt{f}}

The rest of the proof of \rt{f} will be given along the lines of the proof of \rt{h}. We   put
 $$
r_{m+1}=(1-2\delta_m)r_m,\qquad \delta_m=\frac{e^{-2}n^{-1}}{2(m+2)^2},\quad m=0,1,\ldots.
$$
 The following  numerical result can be proved as \rl{4.2}.
\bl{4.2r}
 Let $r_m,\delta_m$ be given as above$.$ Assume that $a_0, a_1, \ldots, a_N$  are non-negative numbers 
satisfying
\be{4.7r}
a_{m+1}\leq
\frac{c_7a_m^2}{r_m\delta_m},\quad 0\leq m<N,
\end{equation}
  and  
\be{4.8r}
  a_0\leq r_0\delta_0^2/c_8,
\end{equation} 
for $c_8=\max\{c_7, 4e^{n+2}\pi\},$
then
\be{4.6r}
a_m\leq r_m\delta_m^2/c_8,\quad m\geq 0.
\end{equation}
\end{lemma}
 
For the proof of \rt{f}, we shall find a sequence of  transformations 
$$
\psi_m\colon z_j^\prime =z_jg_j^{(m)}(z), \quad j=1,\ldots, n,
$$ 
and then  put
\be{am}
(1+a^{(m+1)}( z))\Omega=\psi_m^*((1+a^{(m)}( z))\Omega).
\end{equation}
 Denote
$$
 a_m=\| a^{(m)}\|_{r_m} .
$$
We require that  \re{4.6r} holds and 
$\psi_m$ satisfies
\be{gm}
\|\log g_j^{(m)}\|_{r_{m+1}}\leq r_m\delta_m^2.
\end{equation}
 In particular, \re{gm} implies that
\be{mr}
\psi_m \colon A_{r_{m+1}}\to A_{r_m}.
\end{equation}

 To determine $\psi_m$, we now let $a^{(0)}(z)=a(z) $ be as in \rt{f}.  Then for $\epsilon$ small, we know that $a_0$ satisfies 
\re{4.8r}.  By applying \rl{f1} to $a=a^{(0)}$ and $r=r_0$,  we find a transformation 
$\psi_0=\varphi_{-1}$ as in 
\rl{f1}.
  Next,  we define $a^{(1)}(z)$ by \re{am}. 
From \rp{f2}, we know that $a_0$ and $a_1$ satisfies \re{4.7r}. Hence, \rl{4.2r} implies that \re{4.6r} holds for $a_1$.
 In return, we construct a transformation $\psi_1=\psi$ as in \rl{f1} for $a=a^{(1)}$ and $r=r_1$, so $a^{(2)}(z)$ is determined by \re{am}.  Through this recursive process, we can construct a sequence of transformations $\psi_m$  such that \re{4.6r}-\re{mr} hold.

Notice that $r_m>r_0/2$ for all $m$. Then \re{mr} implies that 
$\psi_0\circ\ldots\circ\psi_m$ is well-defined on the domain $A_{r_0/2}$.   We have
$$
\|\psi_k-\hbox{Id}\|_{r_{k+1}}\leq e\max_{1\leq j\leq n}\{
\|g_j^{(k)}-1\|_{r_{k+1}}\}\leq 
e^2r_k\delta_k^2,
$$
where the last inequality comes from \re{gm} and \re{exp}. Hence, \re{mr} implies that
$$
\|\psi_k\circ\ldots\circ\psi_m-\psi_{k+1}\circ\ldots\circ\psi_m\|_{r_0/2}
\leq \|\psi_k-\hbox{Id}\|_{r_{k+1}}\leq e^2r_0\delta_k^2.
$$
Thus,  the sequence $\psi_0\circ\ldots\circ\psi_m$ converges  on $A_{r_0/2}$, of which the limit mapping  $\psi$ satisfies 
 \be{j1}
\|\psi-\mbox{Id}\|_{r_0/2}\leq e^2r_0\sum_{m=0}^\infty \delta_m^2\leq \frac{r_0}{16n^2e^2},
\end{equation}
where the last inequality is obtained from the choice of $\delta_m$ and the elementary inequality 
$(m+2)^4\geq 4(m+1)(m+2)$. 

It is clear that $\psi^*\omega=\Omega$ on $A_{r_0/2}$. We now can  complete the proof of
\rt{f} by showing that $\psi$ has an inverse mapping $\varphi\colon A_{r_0/8}\to A_{r_0/4}$. To this
end, we fix $z\in  A_{r_0/8}$ and consider the mapping 
$$
T\colon \xi\to z-\psi(\xi)+\xi.
$$
From \re{j1}, it is clear that $T\colon A_{r_0/4}\to A_{r_0/4}$. Next, we want to show that $T$ is
a contraction mapping. For if $\xi^\prime, \xi^{\prime\prime}$ are 
two distinct points in $A_{r_0/4}$, 
we choose $\theta^\prime,\theta^{\prime\prime}\in S_{r_0/4}$ such that $\pi(\theta^\prime)=\xi^\prime$
and $\pi(\theta^{\prime\prime}) =\xi^{\prime\prime}$. We may assume that
$$
|\mbox{Re}\, (\theta_j^{\prime\prime}-\theta_j^\prime)|\leq\pi,\quad j=1,\ldots,n.
$$
Let $\gamma(t)=\pi((1-t)\theta^\prime+t\theta^{\prime\prime})$. Then $\gamma\colon [0,1]\to A_{r_0/4}$, and 
$$
\|\gamma^\prime(t)\|\leq e\|\theta^{\prime\prime}-\theta^\prime\|,\quad t\in[0,1].
$$
Using Cauchy inequalities, we also obtain from \re{j1} that
\be{j3}
\|D\psi-\mbox{Id}\|_{r_0/4}\leq \frac{1}{4en}.
\end{equation}
Combining with the formula
$$
T(\xi^{\prime\prime})-T(\xi^\prime)=-\int_0^1 (D\psi-\mbox{Id})\gamma^\prime(t)\, dt,
$$
we obtain
\be{j2}
\|T(\xi^{\prime\prime})-T(\xi^\prime)\|\leq \|\theta^{\prime\prime}-\theta^\prime\|/4.
\end{equation}
Since $|\xi_j^\prime|, |\xi_j^{\prime\prime}|\geq r_0/4$, then we have
$$
|\xi_j^{\prime\prime}-\xi_j^\prime|\geq \big| |\xi_j^{\prime\prime}|-|\xi_j^\prime|\big|\geq e^{-r_0/4}
|\mbox{Im}\, (\theta_j^{\prime\prime}-\theta_j^\prime)|
$$
and
$$
|\xi_j^{\prime\prime}-\xi_j^\prime|\geq e^{-r_0/4}
\big| \xi_j^{\prime\prime}/|\xi_j^{\prime\prime}|-\xi_j^\prime/|\xi_j^\prime|\big|
\geq 2e^{-r_0/4}|\mbox{Re}\,(\theta_j^{\prime\prime}-\theta_j^\prime)|/\pi.
$$
Thus, we get
$$
\|\xi^{\prime\prime}-\xi^\prime\|\geq \sqrt{2}e^{-1/4}
\|\theta^{\prime\prime}-\theta^\prime\|/\pi.
$$
Combining with \re{j2}, we obtain
$$
\|T(\xi^{\prime\prime})-T(\xi^\prime)\|\leq
\frac{\pi e^{1/4}}{4\sqrt{2}}\|\xi^{\prime\prime}-\xi^\prime\|,\quad
\xi^\prime, \xi^{\prime\prime}\in A_{r_0/4}.
$$
This shows that $T\colon A_{r_0/4}\to A_{r_0/4}$ is a contraction mapping. Therefore, the fixed point theorem implies that $\psi$ has a unique inverse  $\varphi\colon A_{r_0/8}\to A_{r_0/4}$. From \re{j3}, it is also clear that $\varphi$ is holomorphic.
The proof of \rt{f} is complete.

 {\bf Acknowledgment.} The author would like to thank Professor Sidney Webster for the constant encouragement and helpful discussions.

 \bibliographystyle{plain}

\end{document}